\documentclass[12pt]{article} \textwidth=6in \oddsidemargin=0in
\usepackage{amssymb} 
\begin{document}\def\ov{\over} 
\newcommand{\C}[1]{{\cal C}_{#1}} \def\inv{^{-1}}\def\be{\begin{equation}} \def\ep{\varepsilon}
\def\ee{\end{equation}}\def\x{\xi}\def\({\left(} \def\){\right)} \def\iy{\infty}  \def\ld{\ldots} \def\Z{\mathbb Z} \def\cd{\cdots}
\def\P{\mathbb P}\def\bs{\backslash} \def\t{\tau} 
\newcommand{\br}[2]{\left[{#1\atop #2}\right]_\t}   
\def\sp{\vspace{2ex}}  \def\l{\ell} \def\s{\sigma}  
\newcommand{\brr}[2]{\left[{#1\atop #2}\right]} 
\def\DI{\Delta I} \def\mb{\boldmath} \def\l{\ell} 
\renewcommand\S{\mathbb S} \def\sg{{\rm sgn}\,} \def\ds{\displaystyle}
\def\L{\Lambda} \def\X{\raisebox{.3ex}{$\chi$}} \newcommand{\xs}[1]{\x_{\s(#1)}} \newcommand{\xl}[1]{\x_{\l_{#1}}} \def\a{\alpha} \def\b{\beta}

\hfill   April 25, 2010

\begin{center}{\bf \large Formulas for Joint Probabilities for\\
\vskip.7ex
the Asymmetric Simple Exclusion Process}\end{center}

\begin{center}{\large\bf Craig A.~Tracy}\\
{\it Department of Mathematics \\
University of California\\
Davis, CA 95616, USA}\end{center}

\begin{center}{\large \bf Harold Widom}\\
{\it Department of Mathematics\\
University of California\\
Santa Cruz, CA 95064, USA}\end{center}

\begin{abstract}
In earlier work \cite{TW1} the authors obtained integral formulas for
probabilities for a single particle in the asymmetric simple exclusion
process. Here formulas are obtained for joint probabilities for several particles. In the case of a single particle the derivation here is simpler than the one in the earlier work for one of its main results. 
\end{abstract}

\begin{center}{\bf I. Introduction}\end{center}

The \textit{asymmetric simple exclusion process} (ASEP) 
is a continuous time Markov process of interacting particles on the integer lattice $\Z$. Each particle waits exponential time; then with probability $p$ it jumps one step to the right if the site is unoccupied and otherwise it stays put; with probability $q=1-p$ it jumps one step to the left if the site is unoccupied and otherwise it stays put. We refer the reader to Liggett \cite{Li1, Li2} for a precise definition of the model.

When particles can jump only to the right or only to the left, the model is referred to as the 
\textit{totally asymmetric simple exclusion process} (TASEP).
This process is quite special in that it, unlike ASEP, is a \textit{determinantal process} \cite{Bo, Joh1, Sch1, Sos}.  For example, in $N$-particle TASEP  the  transition probability from an
initial state $Y=\{y_1,\ldots, y_N\}$ to a state $X=\{x_1,\ldots, x_N\}$ in time $t$ is given by an $N\times N$ determinant \cite{Sch1}.   For TASEP with step initial condition ($y_i=i$) Johansson, in the seminal paper \cite{Joh1},   proved a  limit law
for the particle current fluctuations.\footnote{The limiting distribution is the distribution function $F_2$ of random  matrix theory.}  

Subsequently, these particle current limit laws for TASEP were extended to other initial conditions including  periodic initial conditions 
\cite{BFPS}, stationary initial conditions \cite{FS1} and the general two-sided Bernoulli initial condition \cite{BAC, FS1}.  
Much recent work on TASEP has focused on the \textit{joint distributions} of the particle current \cite{BFP, BFPa, BFS, CFP, Joh2, Joh3, PS1}.  In particular for step initial condition,
the joint distribution of the associated height function and its KPZ scaling\footnote{KPZ refers to Kardar, Parisi, and Zhang \cite{KPZ}.} to the $\textrm{Airy}_2$ process \cite{Joh2, PS2} is now well developed  for TASEP \cite{Joh3}. We refer the reader to \cite{FS2} for a recent review of TASEP and closely related random growth models.  

For ASEP the results  have been restricted to probabilities for a single particle  \cite{TW1, TW2, TW5, TW3}.  In \cite{TW2, TW5} the present
authors extended Johansson's limit law for TASEP to ASEP thus proving a stronger form of KPZ universality.
Using  results from \cite{TW2}, Sasamoto and Spohn \cite{SS1, SS2, SS3} and Amir, Corwin and Quastel \cite{ACQ} have obtained explicit formulas for
the exact height distribution for the KPZ equation with narrow wedge initial condition.  This
height distribution  interpolates between a standard Gaussian distribution for small time and the $F_2$ distribution for large time.  To extend the results
of \cite{ACQ, SS1,SS2, SS3} and show that the $\textrm{Airy}_2$ process is the large time limit of KPZ with narrow
wedge initial condition we would first try to extend the results of \cite{TW1,TW2} to joint distributions of particle positions in ASEP.  In this paper we take a \textit{first step} in this direction  in that we derive integral formulas for joint probabilities.  Whether these formulas can be used to study the KPZ scaling limit remains to be seen.
 
Here is an outline of the paper. We consider an $N$-particle system with initial configuration $Y=\{y_1,\ld,y_N\}$, where $y_1<\cd< y_N$, and denote by $\{x_1(t),\ld,x_N(t)\}$ the configuration of the system at time $t$. Our starting point is Theorem 2.1 of \cite{TW1}, a formula for the probability
\be\P_Y(x_i(t)=x_i,\ i=1,\ld,N)\label{P1}\ee
(the subscript indicating the initial configuration),
which is a sum of $N!$ multiple integrals, one for each permutation in the symmetric group $\S_N$. The integrals are over small circles about the origin. In Sec.~II we obtain a formula for
\be\P_Y(x_i(t)=x_i,\ i=1,\ld,m),\label{P2}\ee
the joint probability for the first $m$ particles, by summing (\ref{P1}) over all $x_i$ with $i>m$ which satisfy $x_m<x_{m+1}<\cd<x_N<\iy$. We may do this because the integrals are over small circles. Application of combinatorial identity (1.6) of \cite{TW1} replaces the sum of $N!$ integrals by a sum of $N!/(N-m+1)!$ integrals. This is Theorem~1.

To obtain from this a formula for arbitrary consecutive particles,
\be\P_Y(x_i(t)=x_i,\ i=n,\ld,m),\label{P3}\ee
we have to sum over all $x_i$ with $i<n$ which satisfy $-\iy<x_1<\cd<x_{n-1}<x_n$, and to do this requires first replacing integrals over small circles by integrals over large ones. This is accomplished with the help of Lemma~3.1 of \cite{TW1}, and gives Theorem~2. Then summing over the indicated $x_i$ and applying combinatorial identity (1.7) of \cite{TW1} give a formula for (\ref{P3}) as a sum of integrals over large circles. This is Theorem~3.

One could obtain a formula for the general joint probability
\[\P_Y(x_{m_1}(t)=x_{m_1},\ld,x_{m_r}(t)=x_{m_r})\] 
by taking $n=m_1$ and $m=m_r$ in Theorem 3 and summing over those $x_i$ satisfying $x_{m_j}<x_i<x_{m_{j+1}}$. These are finite sums and lead to a rather complicated expression which we shall not write down. There does not seem to be a combinatorial identity that simplifies it.

Finally we consider the special case $\P_Y(x_m(t)=x_m)$,
the probability for a single particle. The formula for consecutive particles simplifies by use of combinatorial identity (1.9) of \cite{TW1}, and yields Theorem~4. This is exactly Theorem~5.2 of \cite{TW1}. The proof here is much simpler than the one in \cite{TW1}, although they use the same ingredients.\footnote{The theorem was the starting point for some of the advances mentioned above, and so it is useful to have this more straightforward proof.}  The awkwardness of that proof made it seem unlikely at first that one could find formulas for joint probabilities without great effort. 

Although we assumed that the initial configuration $Y$ was finite, Theorems~2--4 extend to configurations that are semi-infinite on the right, in particular to $Y=\Z^+$ (step initial condition). This is shown in the last section. 

\begin{center}{\bf II. The first \mb$m$ particles -- small contours}\end{center}

Theorem~2.1 of \cite{TW1} is the formula, valid when $p\ne0$,
\[\P_Y(x_i(t)=x_i,\ i=1,\ld,N)\]
\be=\sum_{\s\in\S_N}\int_{\C{r}}\cdots\int_{\C{r}} A_\s\,\prod_i\x_{\s(i)}^{x_i}\;\prod_i\x_i^{-y_i-1}\;e^{\,\sum_i\ep(\x_i)\,t}\,d\x_1\cd d\x_N,\footnote{In general the product symbol $\prod$ refers only to the first term after it.}\label{PYN}\ee
where  
\be A_\s=\sg\s\ {\prod\limits_{i<j}f(\x_{\s(i)},\,\x_{\s(j)})
\ov \prod\limits_{i<j}f(\x_i,\,\x_j)},\ \ \ f(\x,\,\x')=p+q\x\x'-\x,\ \ \ \ep(\x)=p\x\inv+q\x-1.\label{A}\ee
Here $\C{r}$ is a circle with center zero and radius $r$, which is so small that all nonzero poles of the integrand lie outside $\C{r}$. (All contour integrals are to be given a factor $1/2\pi i$.)

To obtain $\P_Y(x_i(t)=x_i,\ i=1,\ld,m)$ we sum over all $x_{m+1},\ld,x_N$ such that $x_m<x_{m+1}<\cd<x_N<\iy$, which we may do when $r<1$. The result of the summation is that the product $\prod_i\x_{\s(i)}^{x_i}$ in (\ref{PYN}) is replaced by
\[\prod_{i<m}\x_{\s(i)}^{x_i}\ \prod_{i\ge m}\x_{\s(i)}^{x_m}
\ {\x_{\s(m+1)}\,\x_{\s(m+2)}^2\cd\x_{\s(N)}^{N-m}\ov
(1-\x_{\s(m+1)}\,\x_{\s(m+2)}\cd\x_{\s(N)})\cd(1-\x_{\s(N)})}.\]

We take a fixed ordered $(m-1)$-tuple $\L=\{\l_1,\ld,\l_{m-1}\}$ with distinct $\l_i\le N$ and consider first the sum over all permutations $\s$ such that
\[\s(1)=\l_1,\ld,\s(m-1)=\l_{m-1}.\]
(This meaning of $\L$ as an $(m-1)$-tuple will be retained throughout.)

The product $\prod_{i<m}\x_{\s(i)}^{x_i}\;\prod_{i\ge m}\x_{\s(i)}^{x_m}$ becomes
\[\prod_{i<m}\x_{\l_i}^{x_i}\; \prod_{j\not\in\L}\x_j^{x_m},\]
(in the second factor here, and elsewhere below, $\L$ is thought of as a set) and the numerator in the expression for $A_\s$ may be written
\[\prod_{j\ne\l_1}f(\x_{\l_1},\,\x_j)\ \prod_{j\ne\l_1,\,\l_2}f(\x_{\l_2},\,\x_j)\;\cd \prod_{j\ne\l_1,\ld,\l_{m-1}}f(\x_{\l_{m-1}},\,\x_j)\ \prod_{m\le i<j}f(\x_{\s(i)},\,\x_{\s(j)})\] 
\[=\prod_{i\in\L,\,j\not\in\L}f(\x_i,\,\x_j)\ \prod_{i<j}f(\xl{i},\,\xl{j})\ \prod_{m\le i<j}f(\x_{\s(i)},\,\x_{\s(j)}).\]

As for $\sg\s$ we will consider a new guise for $\s$ as its restriction to 
$\L^c=[1,\,N]\bs \L$. (It may be associated in an obvious way with a permutation in $\S_{N-m+1}$, as $\L$ may be associated in an obvious way with a permutation in $\S_{m-1}$.) What we now call $\sg\s$ will be different from the old $\sg\s$, because the number of inversions in the original $\s$ equals the number of inversions of the new guise $\s$ plus the number of inversions in the original $\s$ which involve some $\l_i$. This number is
\[\sum_{i<m}(\l_i-i)+\textrm{\# of inversions in}\ \L.\]
Thus the original $\sg\s$ in (\ref{A}) is to be replaced by
\[(-1)^{\sum_{i<m}(\l_i-i)}\ \sg\L\ \sg\s,\]
where now $\s$ is in its new guise.
 
Identity (1.6) of \cite{TW1} tells us that, with the new guise $\s$,
\[\sum_{\s} \sg\s\,\prod_{m\le i<j}f(\x_{\s(i)},\,\x_{\s(j)})\ 
{\x_{\s(m+1)}\,\x_{\s(m+2)}^2\cd\x_{\s(N)}^{N-m}\ov
(1-\x_{\s(m+1)}\,\x_{\s(m+2)}\cd\x_{\s(N)})\cd(1-\x_{\s(N)})}\]
\[=p^{(N-m+1)(N-m)/2}\;\Big(1-\prod_{j\not\in\L}\x_j\Big) {\ds{\prod_{{i<j\atop i,\,j\not\in\L}}(\x_j-\x_i)}\ov \ds{\prod_{j\not\in\L}(1-\x_j)}}.\]  

Now we sum over all ordered $(m-1)$-tuples $\L$ and obtain

{\bf Theorem 1}. We have when $p\ne0$,\footnote{From now until the last section $N=|Y|$, all indices belong to $[1,\,N]$, and all sets of indices are subsets of $[1,\,N]$.}  
\[\P_Y(x_i(t)=x_i,\ i=1,\ld,m)=\sum_{\L}(-1)^{\sum_{i<m}(\l_i-i)}\ \sg\L\ p^{(N-m+1)(N-m)/2}\]
\[\times\int_{\C{r}}\cdots\int_{\C{r}}
{\ds{\prod_{i\in\L,\,j\not\in\L}f(\x_i,\,\x_j)\ \prod_{i<j}f(\xl{i},\,\xl{j})}\ov
\ds{\prod_{i<j}f(\x_i,\,\x_j)}}\]
\[\times\;\Big(1-\prod_{j\not\in\L}\;\x_j\Big)\; {\ds{\prod_{{i<j\atop i,\,j\not\in\L}}(\x_j-\x_i)}\ov \ds{\prod_{j\not\in\L}(1-\x_j)}}\;\prod_{i< m}\x_{\l_i}^{x_i}\ \prod_{j\not\in\L}\x_j^{x_m}\ \prod_i\x_i^{-y_i-1}\ e^{\,\sum_i\ep(\x_i)\,t}\ d\x_1\cd d\x_N.\]

\begin{center}{\bf III. The first \mb$m$ particles -- large contours}\end{center} 

Here we shall find another representation of the same probability in which all contours of integration are $\C{R}$, where $R$ is arbitrarily large. Observe that some factors $f(\x_i,\,\x_j)$ in the denominator in the integrand are cancelled by factors in the numerator. We assume at first that $p,\,q\ne0$.

{\it Step 1}. We show first that we may take all $\x_\l$-contours with $\l\in\L$ to be $\C{R}$. Let us first expand the $\xl{\a}$-contour, where $\l_\a=\max\L$. There appear to be poles at
\[\xl{\a}={p\ov1-q\x_j}\]
coming from the factor $f(\xl{a},\x_j)$ in the denominator. Since 
$\l_\a=\max\L$ we must have $j\not\in\L$. But then this factor is cancelled by the same factor in the numerator. So the only poles when we expand the $\xl{\a}$-contour are at
\[\xl{\a}={\x_i-p\ov q\x_i},\]
which comes from the factor $f(\x_i,\,\xl{a})$ in the denominator when $i<\l_\a$.

We show that the integral with respect to $\x_i$ of the residue at this pole  equals zero. When we make the substitution $\xl{\a}\to (\x_i-p)/q\x_i$ we find that $\ep(\x_i)+\ep(\xl{\a})$ becomes analytic at $\x_i=0$. As for the first quotient in the integrand, each factor $f$ having $\xl{\a}$ as one of its variables is of the order $\x_i\inv$ as $\x_i\to0$. There are $N-1$ such factors in the numerator and $N-2$ in the denominator (since we don't include the factor $f(\x_i,\,\xl{a})$ that gave rise to the pole). Thus this quotient is $O(\xi_i\inv)$. The residue of $1/f(\x_i,\,\xl{a})$ at the pole is $1/q\x_i$, so these combine to give the power $\x_i^{-2}$. As for the product of powers of the variables at the end of the integrand, when $i\not\in\L$ the product of those involving $\x_i$ is
\[\({\x_i-p\ov q\x_i}\)^{x_\a}\,\x_i^{x_m}\,\({\x_i-p\ov q\x_i}\)^{-y_{\l_\a}}\,\x_i^{-y_i-1}=O(\x_i^{-x_\a+x_m+y_{\l_\a}-y_i})=O(\x_i^2)\]
since $x_\a<x_m$, and $i<\l_\a$. Thus the integrand is analytic inside $\C{r}$ and so its integral is zero. If $i\in\L$ then $i=\ell_{\beta}$ with $\b>\a$. (Otherwise the factor $f(\xl{\b},\,\xl{\a})$ in the denominator is cancelled by the same factor in the numerator.) Then the exponent $x_m$ above is replaced by $x_\b$, and $x_\b>x_\a$.

We have shown that we may expand the $\xl{\a}$-contour to $\C{R}$. Next we expand the $\xl{\b}$-contour where $\l_\b$ is the second-largest element of $\L$. The only difference from what went before is that there could be a pole at
\[\xl{\b}={p\ov1-q\xl{\a}}\]
coming from the factor $f(\xl{\b},\,\xl{\a})$ in the denominator. But $\xl{\a}\in\C{R}$, and if $R$ is large enough the pole would lie inside $\C{r}$ and so would not be passed in the expansion of the $\xl{\b}$-contour. Thus we may expand the $\xl{\b}$-contour to $\C{R}$

Continuing in this way we see that all the $\x_\l$-contours with $\l\in\L$ may be taken to be $\C{R}$.

{\it Step 2}. Next we expand the $\x_i$-contours with $i\not\in\L$, and when we do this we encounter poles at $\x_i=1$ as well as those coming from the first quotient in the integrand. As above the latter poles will not contribute, but those from $\x_i=1$ will. An application of Lemma 3.1 of \cite{TW1} will tell us the result.

We restate the lemma to make it compatible with the present notation. Let $T$ be a finite set of indices and let
\[g=g(\x_i)_{i\in T}\] 
be a function that is analytic for all $\x_i\ne0$. Assume that for $i>k$
\be g\Big|_{\x_i\to (\x_k-p)/q\x_k}=O(\x_k)\label{gcond}\ee
as $\x_k\to0$, uniformly when all the $\x_j$ with $j\ne i,\,k$ are bounded and bounded away from zero. For $S\subset T$ denote by $g_S$ the function obtained from $g$ by setting all $\x_i$ with $i\not\in S$ equal to 1. (In particular $g=g_T$.) Define
\[I_S(\x)=\prod_{{i<j\atop i,j\in S}}{\x_j-\x_i\ov f(\x_i,\x_j)}\,{g_S(\x)\ov\prod_{i\in S}(1-\x_i)}.\]
Then when $p,\,q\ne0$,
\be\int_{\C{r}^{|T|}} I_T(\x)\,\prod_{i\in T}d\x_i=\sum_{S\subset T}
\t^{\s(S,\,T)-|S|}\;{q^{|S|(|S|-1)/2}\ov p^{|T|(|T|-1)/2}}\,
\int_{\C{R}^{|S|}} I_S(\x)\,\prod_{i\in S}d\x_i,\label{gint}\ee
where $R$ is so large that all the zeros of the denominators lie inside $\C{R}$. (When $S$ is empty the integral on the right side is interpreted as $g(1,\ld,1)$.) Here $\t=p/q$ and 
\[\s(S,\,T)=\#\{(i,j): i\in S,\ j\in T,\ i\ge j\}.\]

We apply this with $T=\L^c$ and 
\[g(\x_i)_{i\in\L^c}=\Big(1-\prod_{i\in\L^c}\;\x_i\Big)\ \prod_{i,\,j\in\L^c}\x_i^{x_m-y_i-1}\ e^{\,\sum_{i\in\L^c}\ep(\x_i)\,t}\]
\[\times\;\int_{\C{R}^{|\L|}}{\ds{\prod_{i\in\L,\,j\not\in\L}f(\x_i,\,\x_j)\ \prod_{i<j}f(\xl{i},\,\xl{j})}\ov
\ds{\prod_{{i<j\atop i\,{\rm or}\,j\in\L}}f(\x_i,\,\x_j)}}
\ \prod_{i< m}\x_{\l_i}^{x_i}\ \prod_{i\in\L}\x_i^{-y_i-1}\,e^{\,\sum_{i\in\L}\ep(\x_i)\,t}\,\prod_{i\in\L}d\x_i.\]
The integral on the left side of (\ref{gint}) is then the integral in Theorem 1. That $g$ is analytic for $\x_i\ne0$ despite the denominator in the integrand follows from the fact that $R$ may be arbitrarily large. The argument leading to (\ref{gcond}) is like the analogous argument in Step 1. There are $m-1$ factors of order $\x_k\inv$ in both the numerator and denominator of the quotient in the integrand, and the product $\x_i^{-y_i-1}\,\x_k^{-y_k-1}$ becomes
\[\({\x_k-p\ov q\x_k}\)^{-y_i-1}\,\x_k^{-y_k-1}=O(\x_k)\]
since $i>k$.

So the lemma applies. If $|S|=k$, then the coefficient on the right side of (\ref{gint}) is 
\be p^{-(N-m+1)(N-m)/2}\,\t^{\s(S,\L^c)-k}\,q^{k(k-1)/2}.\label{coeff1}\ee
However there is another factor arising from
\[{\ds{\prod_{i\in\L,\,j\not\in\L}f(\x_i,\,\x_j)}\ov
\ds{\prod_{{i<j\atop i\,{\rm or}\,j\in\L}}f(\x_i,\,\x_j)}}\]
when any $\x_k$ with $k\in\L^c\bs S$ is set equal to one. This factor is
\[{\prod_{i\in\L}f(\x_i,1)\ov \ds{\prod_{{k<j\atop j\in\L}}f(1,\x_j)\cdot\prod_{{i<k\atop i\in\L}}f(\x_i,1)}}=\ds{\prod_{{j>k\atop j\in\L}}{f(\x_j,1)\ov f(1,\x_j)}}=(-\t)^{\#\{j\in\L:\,j>k\}}.\]
Hence the total factor when the $\x_k$ with $k\in\L^c\bs S$ are set equal to one is
\be(-\t)^{\s(\L,\,S^c\cap\L^c)}.\label{coeff2}\ee
If we use the bilinearity of $\s(U,\,V)$ in its two variables, and the easy fact
\[\s(U,\,V)+\s(V,\,U)=|U|\,|V|+|U\cap V|,\]
we see that 
\[\s(\L,\,S^c\cap\L^c)=\s(\L,\,\L^c)-\s(\L,\,S)=\sum\l_i-m(m-1)/2+
\s(S,\,\L)-(m-1)k.\]
Thus (\ref{coeff1}) times (\ref{coeff2}) equals
\[p^{-(N-m+1)(N-m)/2}\,\,(-1)^{\sum_{i<m}(\l_i-i)+\s(S,\,\L)-(m-1)k}\,
\t^{\sum_{i<m}(\l_i-i)+\s(S)-mk}\,q^{k(k-1)/2},\]
where we have written $\s(S)$ for $\s(S,\,\L)+\s(S,\,\L^c)=\s(S,\,[1,\,N])$.

We replace each integral in the formula of Theorem~1 by the sum over $S\subset\L^c$ we  obtain using (\ref{gint}). The result is
\pagebreak 

{\bf Theorem 2}. We have when $q\ne0$, 
\[\P_Y(x_i(t)=x_i,\ i=1,\ld,m)=\sum_{{\L\atop S\subset\L^c}}(-1)^{\s(S,\,\L)-(m-1)k}\ \sg\L\ \t^{\sum_{i<m}(\l_i-i)+\s(S)-mk}\,q^{k(k-1)/2}\]
\[\times\int_{\C{R}}\cdots\int_{\C{R}}
{\ds{\prod_{i\in\L,\,j\in S}f(\x_i,\,\x_j)\ \prod_{i<j}f(\xl{i},\,\xl{j})}\ov
\ds{\prod_{i<j}f(\x_i,\,\x_j)}}\]
\[\times\;\Big(1-\prod_{j\in S}\;\x_j\Big) {\ds{\prod_{{i<j\atop i,\,j\in S}}(\x_j-\x_i)}\ov \ds{\prod_{j\in S}(1-\x_j)}}
\prod_{i<m}\x_{\l_i}^{x_i}\ \prod_{j\in S}\x_j^{x_m}\ \prod_i\x_i^{-y_i-1}\,e^{\,\sum_i\ep(\x_i)\,t}\,\prod_{i}d\x_{i}.\]
The summation is over all sets $S$ and $(m-1)$-tuples $\L$ disjoint from $S$. Here $k=|S|$ and in the integrand indices that are not specified belong to $S\cup\L$.

We assumed at first that $p,\,q\ne0$ since that was required by Lemma 3.1 of \cite{TW1}. To obtain the relation when $p=0$ by passing to the $p \to0$ limit we need only observe that the power of $\t$ in the coefficient is nonnegative and that the integrand is continuous at $p=0$ when $R>1$.

\begin{center}{\bf IV. Consecutive particles}\end{center}

Here we find a formula for $\P_Y(x_i(t)=x_i,\ i=n,\ld,m)$. If we take the preceding formula and sum over all $x_1,\ld,x_{n-1}$ with $-\iy<x_1<\cd<x_{n-1}<x_n$, which we may do when $R>1$, the partial product $\prod_{i<n}\x_{\l_i}^{x_i}$ is replaced by
\be{(\x_{\l_1}\cd\x_{\l_{n-1}})^{x_n}\ov
(\x_{\l_1}-1)(\x_{\l_1}\x_{\l_2}-1)\cd(\x_{\l_1}\cd\x_{\l_{n-1}}-1)}.\label{factor}\ee

Our $(m-1)$-tuple $\L$ may be written in the obvious way as $\{\L_1,\,\L_2\}$, where
\be\L_1=\{\l_1,\ld,\l_{n-1}\},\ \ \ \L_2=\{\l_n,\ld,\l_{m-1}\}.\label{LL}\ee

Until now we did not differentiate in the notation between an ordered set and the corresponding unordered set. Which was meant was always clear from the context. Now we introduce a notation when we want to distinguish the two: for an ordered set $\L$ we denote by $\X(\L)$ the corresponding unordered set.

In the integrand in Theorem 2 the only factor that depends on $\L$ and not just on the set $\X(\L)$ is
\[\prod_{i<j}f(\xl{i},\,\xl{j})=\prod_{i<j<n}f(\xl{i},\,\xl{j})
\ \prod_{i<n,\,j\ge n}f(\xl{i},\,\xl{j})\ \prod_{n\le i<j}f(\xl{i},\,\xl{j}).\]
The first factor depends on $\L_1$, the second depends only on $\X(\L_1)$ and $\X(\L_2)$, while the third depends on $\L_2$.

In Theorem 2 we are to sum over all sets $S$ and all
$(m-1)$-tuples $\L$ disjoint from $S$. Let us fix $S$ first and sum over all $\L$ disjoint from $S$. This sum over $\L$ is the sum over all $\L_1$ and $\L_2$ as in (\ref{LL}). We do this sum by first taking a fixed set $S_1$ disjoint from $S$ with $|S_1|=n-1$, and summing over all $\L_1$ with $\X(\L_1)=S_1$. Then we must sum over $S_1$ and, of course, $\L_2$ (which is disjoint from $S_1$ and $S$) and $S$. To recapitulate:  the summation in the statement of Theorem 2 may be replaced by
\[\sum_{S,\,S_1,\,\L_1,\,\L_2}\]
where $S$ and $S_1$ run over all disjoint sets with $|S_1|=n-1$, and $\L_1$ and $\L_2$ run over all $(n-1)$-tuples and $(m-n)$-tuples respectively satisfying 
\[ \X(\L_1)=S_1,\ \ \ \X(\L_2)\subset(S\cup S_1)^c.\]

With this notation the coefficient in Theorem~2 may be written
\[(-1)^{\s(S\cup S_1,\,\L_2)-(m-1)k}\ (-1)^{\s(S_1,\,\L_2)}\,\sg\L_1\ \sg\L_2\ \t^{\s(S\cup S_1\cup\L_2)-m(m-1)/2-mk}\,q^{k(k-1)/2}.\]
Notice that $\L_1$ appears here only in the factor $\sg\L_1$. 
The sum of those terms involving $\L_1$ as distinguished from $S_1$ is
\[\sum_{\L_1}\sg\,\L_1\,\prod_{i<j<n}f(\xl{i},\,\xl{j})\ {(\x_{\l_1}\cd\x_{\l_{n-1}})^{x_n}\ov
(\x_{\l_1}-1)(\x_{\l_1}\x_{\l_2}-1)\cd(\x_{\l_1}\cd\x_{\l_{n-1}}-1)}.\]
(Recall the factor (\ref{factor}).)
Identity (1.7) of \cite{TW1} tells us that this sum equals $(\x_{\l_1}\cd\x_{\l_{n-1}})^{x_n}$ times
\[q^{(n-1)(n-2)/2}\ {\ds{\prod_{{i<j\atop i,\,j\in S_1}}(\x_j-\x_i)}\ov \ds{\prod_{j\in S_1}(\x_j-1)}}=(-1)^{n-1}\ q^{(n-1)(n-2)/2}\ {\ds{\prod_{{i<j\atop i,\,j\in S_1}}(\x_j-\x_i)}\ov \ds{\prod_{j\in S_1}(1-\x_j)}}.\]

We write the last product here times the analogous product in the integrand in Theorem 2, where $S$ appears instead of $S_1$, as
\[{\ds{\prod_{{i<j\atop i,\,j\in S\cup S_1}}(\x_j-\x_i)}\ov \ds{\prod_{j\in S\cup S_1}(1-\x_j)}}\ {(-1)^{\s(S_1,\,S)}\ov\ds{\prod_{i\in S_1,\ j\in S}(\x_j-\x_i)}}.\]

Putting all these things together and using $\s(S,\,S_1)+\s(S_1,\,S)=(m-1)k$ we obtain

{\bf Theorem 3}. We have when $q\ne0$,
\[\P_Y(x_i(t)=x_i,\ i=n,\ld,m)\]
\[=\sum_{S,\,S_1,\,\L_2}(-1)^{\s(S\cup S_1,\,\L_2)+n-1}\ \sg\L_2\ \t^{\s(S\cup S_1\cup\L_2)-m(m-1)/2-mk}\,q^{k(k-1)/2+(n-1)(n-2)/2}\]
\[\times\int_{\C{R}}\cdots\int_{\C{R}}
{\prod\limits_{i\in S_1,\;j\in S\cup\L_2}f(\x_i,\,\x_j)\ 
\prod\limits_{i\in\L_2,\;j\in S}f(\x_i,\,\x_j)\ \prod\limits_{n\le i<j}f(\x_{\l_i},\,\x_{\l_j})\ov\prod\limits_{i<j}f(\x_i,\,\x_j)}\]
\[\times\Big(1-\prod_{j\in S}\;\x_j\Big) {\ds{\prod_{{i<j\atop i,\,j\in S\cup S_1}}(\x_j-\x_i)}\ov \ds{\prod_{j\in S\cup S_1}(1-\x_j)}}\ {1\ov\ds{\prod_{i\in S_1,\ j\in S}(\x_j-\x_i)}}\]
\[\times 
\prod_{i\in S_1}\x_{i}^{x_n}\ \prod_{n\le i<m}\x_{\l_i}^{x_i}\ \prod_{j\in S}\x_j^{x_m}\ \prod_i\x_i^{-y_i-1}\,e^{\,\sum_i\ep(\x_i)\,t}\,\prod_{i}d\x_{i}.\]
The summation is over all disjoint sets $S$ and $S_1$ with $|S_1|=n-1$ and $(m-n)$-tuples $\L_2$ disjoint from $S\cup S_1$. Here $k=|S|$ and 
in the integrand indices that are not specified belong to $S\cup S_1\cup\L_2$.

\begin{center}{\bf V. One particle}\end{center}

In this case $n=m$ and so $\L_2$ disappears, the first line of the right side becomes
\[\sum_{S,\,S_1}(-1)^{m-1}\ \t^{\s(S\cup S_1)-m(m-1)/2-mk}\,q^{k(k-1)/2+(m-1)(m-2)/2},\]
and the factor
\[\prod_{i\in S_1}\x_{i}^{x_n}\ \prod_{n\le i<m}\x_{\l_i}^{x_i}\ \prod_{j\in S}\x_j^{x_m}\]
becomes
\[\prod_{i\in S\cup S_1}\x_{i}^{x_m}.\]

Now we take a fixed set $S_3$ with $|S_3|=m+k-1$ and sum over all partitions $S_3=\{S,\,S_1\}$ with $|S|=k,\ |S_1|=m-1$. The only part of the sum that depends on $S$ and $S_1$ individually is
\[\sum_{S,\,S_1}\ {\prod\limits_{i\in S_1,\,j\in S}f(\x_i,\,\x_j)\ov\ds{\prod_{i\in S_1,\ j\in S}(\x_j-\x_i)}}\ \Big(1-\prod_{j\in S}\;\x_j\Big).\]
If we observe that $S$ is the complement of $S_1$ in $S_3$ we see that we can apply identity (1.9) of \cite{TW1} (with $S$ there $S_1$ here, with $m$ there $m-1$ here, and with $N$ there $m+k-1$ here), which tells us that the sum equals
\[q^{(m-1)k}\,\br{m+k-2}{m-1}\Big(1-\prod_{i\in S_3}\x_i\Big),\]
where the $\t$-binomial coefficient $\br{N}{n}$ is defined by
\[\br{N}{n}={(1-\t^N)\,(1-\t^{N-1})\cdots (1-\t^{N-n+1})\ov (1-\t)\,(1-\t^2)\cdots (1-\t^n)}.\]
Hence, after some algebra, our formula becomes

{\bf Theorem 4}. We have when $q\ne0$,
\[\P_Y(x_m(t)=x_m)=(-1)^{m-1}\,\t^{m(m-1)/2}\,\sum_{|S_3|\ge m}\,\t^{\s(S_3)-m|S_3|}\,q^{|S_3|(|S_3|-1)/2}\,\br{|S_3|-1}{m-1}\]
\[\times \int_{\C{R}}\cdots\int_{\C{R}}\ 
\prod_{i<j}{\x_j-\x_i\ov f(\x_i,\,\x_j)}\ {1-\prod_i\x_i\ov 
\prod_i(1-\x_i)}\
\prod_i\x_{i}^{x_m-y_i-1}\,e^{\,\sum_i\ep(\x_i)\,t}\,\prod_{i}d\x_{i},\]
where all indices in the integrand run over $S_3$.

This is exactly Theorem 5.2 of \cite{TW1}.

\begin{center}{\bf VI. Semi-infinite configurations}\end{center}

Here we show that Theorem 2, and consequently also Theorems 3 and 4, hold when $Y$ is semi-infinite on the right, if we take the sums
over finite sets $S\subset\L^c$. It follows from the fact that ASEP is a Feller process \cite{Li1} that the probability for the semi-infinite initial configuration $Y$ equals the $N\to\iy$ limit of the probability for initial configuration $\{y_1,\ld,y_N\}$. Thus we sum over only those $\L\subset[1,\,N]$ and those $S\subset[1,\,N]$ in Theorem~2 and then pass to the limit. The limit will be the sum over all $\L$ and (finite) $S$ if the resulting series is absolutely convergent. This we now show, using the fact that $R$ may be taken arbitrarily large. 

Consider the integrand first. Each $f$ factor is of the order $R^2$ for large $R$ while each $\x_i$ is of order $R$. Combining the estimates for all factors in the integrand other than $\prod_i\x_i^{-y_i-1}$ gives, after a little algebra, $O(R^{-k^2/2+O(k)})$. (Recall that $m$ is fixed.) Since each $y_i\ge y_1+i-1$, we have $\prod_i\x_i^{-y_i-1}=O(R^{-\s(\L\cup S)+O(k)})$. Thus the integrand is
\[O(R^{-k^2/2-\s(\L\cup S)+O(k)}).\]
Another factor $R^k$ comes from the domain of integration, but this does not change the bound.

The coefficient on the right side of Theorem 2 has absolute value at most $\t^{\s(\L\cup S)+O(k)}$. If we take $R>\t^2$, as we may, then this combined with the preceding bound gives
\be O(R^{-k^2/2-\s(\L\cup S)/2+O(k)}).\label{bound}\ee

This a bound for the summand in Theorem 2. Now we show that the sum of these over all $\L$ and $S$ is finite. For any $s$ we have
\[\#\{(\L,\,S):\s(\L\cup S)=s\}\le 2^s\,s^{m-1}.\]
The reason is that if $\s(\L\cup S)=s$ then the largest element in $\L\cup S$ is at most $s$. So there are at most $2^s$ choices for $S$, snd having chosen $S$ there are at most $s^{m-1}$ choices for $\L$. It follows that the sum of (\ref{bound}) over all $(\L,\,S)$ is at most a constant times
\[\sum_{k,\,s\ge 0}2^s\,s^{m-1}\,R^{-k^2/2-s/2+O(k)},\]
which is finite when $R>4$.

\begin{center}{\bf Acknowledgment}\end{center}

This work was supported by the National Science Foundation through grants DMS-0906387 (first author) and DMS-0854934 (second author).

\end{document}